\documentclass[12pt,reqno]{amsart}
\usepackage{xy, amssymb, amsfonts, amsbsy, amsthm, amsmath, amscd, latexsym, stmaryrd, epic, eepic,eucal}

\DeclareMathAlphabet{\mathpzc}{OT1}{pzc}{m}{it}

\newenvironment{dem}{\begin{proof}[\bf Proof]}{\end{proof}}

\newtheorem{theorem}{\bf Theorem}[section]
\newtheorem{lemma}[theorem]{\bf Lemma}
\newtheorem{propos}[theorem]{\bf Proposition}

\newtheorem{claim}[theorem]{\bf Claim}

\theoremstyle{definition}

\newtheorem{oss}[theorem]{\bf Remark}

\newcommand{\A}{\mathbb A}
\newcommand{\B}{\text{B}}

\newcommand{\Cg}{\mathfrak C}

\newcommand{\Gm}{\mathbb G_{\textbf{m}}}

\newcommand{\M}{\mathfrak M}

\newcommand{\Ms}{\mathcal M}

\newcommand{\Nf}{\mathcal N}

\newcommand{\Pro}{\mathbb P}
\newcommand{\Q}{\mathbb Q}

\newcommand{\aut}{\text{Aut}}

\newcommand{\cart}{\ar @{} [dr] |{\Box}}

\newcommand{\km}{\mathpzc k}

\newcommand{\rat}{\otimes \mathbb Q}

\newcommand{\sla}{\mathfrak{sl}}

\newcommand{\unoa}{\xymatrix { *=0{\bullet} \ar@{-}[r] & *=0{\bullet} }}
\newcommand{\duea}{\xymatrix { *=0{\bullet} \ar@{-}[r] & *=0{\bullet} \ar@{-}[r] &*=0{\bullet} }}

\newcommand{\trea}{\xymatrix { *=0{\bullet} \ar@{-}[r] & *=0{\bullet} \ar@{-}[r] &*=0{\bullet} \ar@{-}[r] &*=0{\bullet} }}

\newcommand{\treb} {\xymatrix @R=8pt { & *=0{\bullet} \ar@{-}[dd] & \\ & &  \\ & *=0{\bullet} \ar@{-}[dl] \ar@{-}[dr] & \\ *=0{\bullet} & & *=0{\bullet}}}

\newcommand{\trebg}{\xymatrix @R=8pt { & *=0{} \ar@{-}[4,0] & \\ *=0{} \ar@{-}[0,2] & & *=0{} \\ *=0{} \ar@{-}[0,2] & & *=0{}\\ *=0{} \ar@{-}[0,2] & & *=0{}\\ & *=0{} &}}

\newcommand{\quattroa}{\xymatrix { *=0{\bullet} \ar@{-}[r] & *=0{\bullet} \ar@{-}[r] &*=0{\bullet} \ar@{-}[r] &*=0{\bullet} \ar@{-}[r] &*=0{\bullet} }}

\newcommand{\quattrob}{\xymatrix {& *=0{\bullet}  \ar@{-}[d]& \\ *=0{\bullet} \ar@{-}[r] & *=0{\bullet} \ar@{-}[r] \ar@{-}[d] & *=0{\bullet}\\ & *=0{\bullet} &}}

\newcommand{\quattroc}{\xymatrix @R=4pt {*=0{\bullet} \ar@{-}[dr] & & &\\ & *=0{\bullet} \ar@{-}[r] &*=0{\bullet} \ar@{-}[r] &*=0{\bullet}\\ *=0{\bullet} \ar@{-}[ur] & & &}}

\newcommand{\quattrocg}{\xymatrix @R=8pt { & *=0{} \ar@{-}[4,0] & & &\\ *=0{} \ar@{-}[0,2] & & *=0{} &*=0{} \ar@{-}[2,0] & \\ *=0{} \ar@{-}[0,4] & & & &*=0{}\\ *=0{} \ar@{-}[0,2] & & *=0{} &*=0{} &\\ & *=0{} & & &}}

\newcommand{\cinquenz}{\xymatrix { *=0{\bullet} \ar@{-}[d] & *=0{\bullet} \ar@{-}[d] \\ *=0{\bullet} \ar@{-}[r] & *=0{\bullet} \\
*=0{\bullet} \ar@{-}[u] & *=0{\bullet} \ar@{-}[u]}}

\input xy
\xyoption{all}

\begin{document}

\title[Chow Ring of Rational Curves]{THE CHOW RING OF THE STACK OF RATIONAL CURVES WITH AT MOST 3 NODES}
\author{Damiano Fulghesu}
\address{Department of Mathematics,
University of Missouri,
Columbia, MO 65211}
\email{damiano@math.missouri.edu} 

\begin{abstract}
In this paper we explicit the rational Chow ring of the stack $\M_{0}^{\leq 3}$ consisting of nodal curves of genus $0$ with at most $3$ nodes: it is a $\Q$-algebra with $10$ generators and $11$ relations.
\end{abstract}

\maketitle

\section{Introduction}

\medskip

In this paper we explicit the ring $A^*(\M_{0}^{\leq 3})\otimes \Q$ (Theorem \ref{finthr}). Our technique is to compute $A^*(\M_{0}^{\leq n})\rat $ for $n \leq 3$ by induction on $n$. For $n = 0$ we have $\M_{0}^{0} = \B \Pro GL_{2}$, and this case is well understood \cite{Pa}.

The inductive step is based on the following fact: if $n \leq 4$, then the top Chern class of the normal bundle of $\M_{0}^{n}$ into $\M_{0}^{\leq n}$ is not a $0$-divisor in $A^*(\M_{0}^{n})\rat$. As a consequence, by an elementary algebraic Lemma (\ref{incollamento}) we can reconstruct the ring $A^*(\M_{0}^{\leq n})\rat $ from the rings $A^*(\M_{0}^{\leq n-1})\rat$ and $A^*(\M_{0}^{n})\rat$, provided that we have an explicit way of extending each class in $A^*(\M_{0}^{\leq n-1})\rat$ to a class in $A^*(\M_{0}^{\leq n})\rat$ (as we have seen in \cite{fulg2}), and then computing the restriction of this extension to $A^*(\M_{0}^{n})\rat$.

We find $10$ generators: the classes $\gamma_{\Gamma}$ for all trees $\Gamma$ with at most three nodes (they are $5$), plus the Mumford class $\km_{2}$. The remaining $4$ generators are somewhat new with respect to the tautological classes introduced for stable curves.

The ideal of relations is determined essentially through two useful technical algebraic Lemmas (\ref{incollamento} and \ref{quoz}). 

\section{Algebraic lemmas} Here we state and prove two algebraic lemmas for future reference. The first one is Lemma 4.4 \cite{VV}:
\begin{lemma} \label{incollamento}
Let $A, B$ and $C$ be rings, $f:B \to A$ and $g:B \to C$ ring homomorphism. Let us suppose that there exists an homomorphism of abelian groups $\phi: A \to B$ such that the sequence
$$
A \xrightarrow{\phi} B \xrightarrow{g} C \xrightarrow{} 0
$$
is exact; the composition $f \circ \phi: A \to A$ is multiplication by a central element $a \in A$ which is not a 0-divisor.

Then $f$ and $g$ induce an isomorphism of ring
$$
(f,g):B \to A \times_{A/(a)} C,
$$
where the homomorphism $A \xrightarrow{p} A/(a)$ is the projection, while $C \xrightarrow{q} A/(a)$ is induced by the isomorphism $C \simeq B/\ker g$ and the homomorphism of rings $f: B \to A$.
\end{lemma}

\begin{dem}
Owing to the fact that $a$ is not a 0-divisor we immediately have that $\phi$ is injective. Let us observe that the map $(f,g):B \to A \times_{A/(a)} C$ is well defined for universal property and for commutation of the diagram:
\begin{equation*}
\xymatrix{
B  \ar[r]^g \ar[d]_{f} &C \ar[d]^{q}\\
A \ar[r]^{p} &A/(a)
}
\end{equation*}
Now let us exhibit the inverse function $A \times_{A/(a)} C \xrightarrow{\rho} B$. Given $(\alpha, \gamma) \in  A \times_{A
/(a)} C$, let us chose an element $\beta \in B$ such that $g(\beta)=\gamma$. By definition of $q$ we have that $f(\beta)-\alpha$ lives in the ideal $(a)$ and so (it is an hypothesis  on $f \circ \phi$) there exists in $A$ an element $\tilde \alpha$ such that:
$$
f(\beta)-\alpha=f(\phi(\tilde \alpha))
$$
from which:
$$
f(\beta - \phi( \tilde \alpha)) = \alpha
$$
we define then $\rho (\alpha,\gamma):=\beta-\phi( \tilde \alpha)$. In order to verify that it is a good definition, let us suppose that there exist an element $\beta_0 \in B$  such that  $(f,g)(\beta_0)=0$, to be precise there exists an element $\alpha_0 \in A$ such that: $\phi(\alpha_0)=0$ and furthermore
$$
0=f(\phi(\alpha_0))=a\alpha_0
$$
but we have that $a$ is not a divisor by zero, so necessary we have $\alpha_0=0$ and $\beta_0=0$. we conclude by noting that from the definition of $\rho$ we have immediately that it is an isomorphism. 
\end{dem}

\begin{oss}
The Lemma \ref{incollamento} will be used for computing the Chow ring of $\M_0^{\leq k+1}$ when the rings $A^*[\M_0^{\leq k}]\rat$ and $A^*[\M_0^{k+1}] \rat$ are known. As a matter of fact, given an Artin stack $\mathcal X$ and a closed Artin substack $\mathcal Y \xrightarrow{i} \mathcal X$ of positive codimension, we have the exact sequence of groups (see \cite{kre} Section 4)
$$
A^*(\mathcal Y) \rat \xrightarrow{i_*} A^*(\mathcal X) \rat \xrightarrow{j^*} A^*(\mathcal U) \xrightarrow{} 0
$$
By using the Self-intersection Formula, it follows that: 

$$
i^*i_*1=i^*[\mathcal Y]=c_{top}(\Nf_{\mathcal Y / \mathcal X}),
$$
when $c_{top}(\Nf_{\mathcal Y / \mathcal X})$ is not 0-divisor we can apply the Lemma.
\end{oss}

We will also use the following algebraic Lemma:

\begin{lemma}\label{quoz}
Given the morphisms

$$
A_1 \xrightarrow{p_1} \overline{A}_1 \xrightarrow{} B \xleftarrow{} \overline{A}_2 \xleftarrow{p_2} A_2
$$
in the category of rings, where the maps $p_1$ and $p_2$ are quotient respectively for ideals  $I_1$ and $I_2$. Then it defines an isomorphism
$$
\overline{A}_1 \times_B \overline{A}_2 \cong \frac{A_1 \times_B A_2}{(I_1,I_2)}.
$$
\end{lemma}
\begin{dem}
Let us consider the map
$$
(p_1, p_2): A_1 \times_B A_2 \to \overline{A}_1 \times_B \overline{A}_2,
$$
by surjectivity of $p_1$ e $p_2$ this map is surjective, while the kernel is the ideal $(I_1,I_2)$.
\end{dem}

\section{The open substack $\M _0 ^0$}

\medskip

\framebox{$\Gamma:= \bullet$}

\medskip

As we have seen in \cite{fulg1}
\begin{equation}
\M_0 ^{\Gamma} \simeq \B \aut(C).
\end{equation}

we have that the stack $\M_0 ^0$ is the classifying space of $\Pro Gl_2$. Owing to the fact that $\Pro Gl_2 \cong SO_3$ and following \cite{Pa} we have
$$
A^*(\M_0 ^0) \rat  \cong \Q[c_2(\sla _2)].
$$
Since (as shown in \cite{fulg2}) $c_2(\sla _2)=(1/2)\km_2$, we can write
\begin{propos}
$$
A^*(\M_0 ^0) \rat =\Q[\km_2].
$$
\end{propos}

\section{The first stratum}

\framebox{$\Gamma:=$ \small{$\unoa$}} 

\medskip

We order the two components. The automorphism group is $\Cg_2 \ltimes (E \times E)$, (where $\Cg_2$ is the order two multiplicative group) and the action of its generator $\tau$ over  $E \times E$ exchanges the components.

Then the induced action of $\tau$ on the ring $A^*_{E \times E} \rat \cong \Q [t_1, t_2]$ exchanges the first Chern classes $t_1$ and $t_2$. The invariant polynomials are the symmetric ones which are algebrically generated by: $\{ (t_1 + t_2)/2, (t_1^2  + t_2 ^2)/2 \}$. By recalling the description of Mumford classes given in \cite{fulg2}, we have $A^*(\M _0 ^1) \rat = \Q[\km_1, \km_2]$.
Let us consider the two inclusions $i$ and $j$
(respectively closed and open immersions) and the \'etale covering $\phi$
\begin{equation}
\xymatrix{
\widetilde \M_0^1 \ar[r]^{\phi}&\M_0 ^1 \ar [dr]^i \\
&\M_0 ^0 \ar [r]_j & \M_0 ^{\leq 1}
}
\end{equation}
we obtain the following exact sequence
$$
A^*(\M_0 ^1) \rat \xrightarrow{i_*} A^*(\M_0 ^{\leq 1}) \rat \xrightarrow{j^*} A^*(\M_0 ^0) \rat \xrightarrow{} 0
$$
for what we have seen we have:
$$
\Q[\km_1, \km_2] \xrightarrow{i_*} A^*(\M_0 ^{\leq 1}) \rat \xrightarrow{j^*} \Q[\km_2] \xrightarrow{} 0.
$$
Now we have that the first Chern class of the normal bundle $N_{\M_0^1} ( \M_0^{\leq 1})$ is
$$
i^*i_* [\M_0^1]=\phi_*\frac{1}{2}(t_1+t_2)=-\km_1
$$
Since $A^*{\M^1_0}$ is an integral domain we can apply Lemma (\ref{incollamento}) and obtain the ring isomomorphism
$$
A^*(\M_0 ^{\leq 1}) \rat \cong \Q[\km_1,\km_2] \times_{\Q[\km_2]} \Q[\km_2] \cong \Q[\km_1, \km_2].
$$
where the map $q: \Q[\km_2] \to \Q[\km_1,\km_2]/(\km_1)=\Q[\km_2]$ tautologically sends $\km_2$ into $\km_2$.

So we have
\begin{propos}
$$
A^*(\M_0 ^{\leq 1}) \rat \cong \Q[\km_1, \km_2].
$$
\end{propos}
\section{The second stratum}

\framebox{$\Gamma:=$ \small{$\duea$}} We order the two components with one node.

In this case the group of automorphism of the fiber is
$$
\text{Aut}(C^{\Gamma}) \cong \Cg_2 \ltimes (\Gm \times E \times E)=:\Cg_2 \ltimes H,
$$
where the action of $\tau$ sends an element $g \in \Gm$ into $g^{-1}$ and exchange the components isomorphic to $E$.

We can identify $A^*(\B (\Gm)^3)$ with $A^*(\widetilde{\M}_0^{\Gamma})$ and $A^*(\B \aut(\Gamma) \ltimes (\Gm)^3)$ with $A^*(\M_0^{\Gamma})$.  Set
\begin{eqnarray*}
t_1=\psi^1(\infty,1) \quad t_2=\psi^1(\infty,2) \quad r=\psi^2(\infty)&
\end{eqnarray*}
the action induced by $\tau$ on these classes is $\tau(r, t_1, t_2) = (-r, t_2, t_1)$.
With reference to the map
$$
\B (\Gm)^3 \xrightarrow{\phi} \B \aut(\Gamma) \ltimes (\Gm)^3
$$
we recall that $\phi^*$ is a ring isomorphism between $A^*(\M_0^1)\rat$ and $A^*(\B (\Gm)^3)^{\Cg_2}$.

We can describe $A^*(\B H) \rat=\Q [r, t_1, t_2]$ as the polynomial ring in $r$ with coefficients in $\Q[t_1,t_2]$, so we write a polynomial  $P(r,t_1,t_2)$ as $\sum _{i=0} ^k r^i P_i(t_1, t_2)$.

The polynomial $P$ is invariant for the action of $\tau$ if and only if the coefficients of the powers of $r$  in $P(r,t_1,t_2)$ are equal to those of the polynomial $P(-r,t_2,t_1)$.

That is to say that $P_i$ with even index are invariant for the exchange of $t_1$ and $t_2$, while those with odd index are anti-invariant. An anti-invariant polynomial $Q$ is such that $Q(t_1,t_2) + Q(t_2, t_1)=0$
and consequently it is the product of $(t_1 - t_2)$ by an invariant polynomial. It is furthermore straightforward verifying that any such polynomial is invariant for the action of $\tau$.

So an algebraic system of generators for $(A^*_{(\Gm)^3}\rat)^{\Cg _2}$ is given by
$$
u_1:=t_1 + t_2 \quad u_2:=t_1 ^ 2 +  t_2 ^2, \quad u_3:=r(t_1-t_2), \quad u_4:=r^2.
$$
We know that
\begin{eqnarray*}
&&\phi^*\km_1=- u_1\\
&&\phi^*\km_2=- u_2\\
&&\phi^*(\gamma_2) = (t_1-r)(t_2+r)=\frac{1}{2}(u_1^2-u_2) +u_3 - u_4
\end{eqnarray*}
where $\gamma_2=c_2(\Nf_{\M_0^2/\M_0^{\leq 2}})$, and there exists a class $x \in A^2\M_0^2\rat$ such that $u_3=\pi^*x$.
\begin{claim}
The ideal of relations is generated on $\Q [\km_1, \km_2, \gamma_2, x]$ by the polynomial
\begin{eqnarray} \label{p2a}
(2x + (2\km_2+\km_1^2))^2-(2\km_2 + \km_1^2)(4\gamma_2 - \km_1^2)=0.
\end{eqnarray}
\end{claim}

\begin{dem}
From direct computation we have that relation (\ref{p2a}) holds and the polynomial is irreducible. On the other hand let us consider the map $f: \A^3_{\Q} \to \A^4_{\Q}$ defined as $(t_1,t_2,r) \mapsto (u_1, u_2, u_3, u_4)$.
If the generic fiber of $f$ is finite then $f(\A^3_{\Q})$ is an hypersurface in $\A^4_{\Q}$ and we have done.
Now for semicontinuity it is sufficient to show that a fiber is finite.
Let us consider the fiber on 0. We have that $u_1, u_2, u_3, u_4$ are simultaneously zero iff $t_1=t_2= r=0$.
\end{dem}

{\bf NOTE:} In the following we do not explicit the argument above.

\bigskip

Now set $\eta:=2x + (2k_2 + k_1^2)$,
we have that $A^*(\M _0 ^2) \rat$ is isomorphic to the graded ring $\Q [\km_1,\km_2,\gamma_2, \eta]/I$
where the ideal $I$ is generated by the polynomial $\eta^2-(2\km_2 + \km_1^2)(4\gamma_2 - \km_1^2)$.
Since $\phi^* \phi_*$ is multiplication by two, we also have the following relation
$$
\eta=\phi_*\left( \frac{1}{2}(t_1-t_2)(2r-t_1+t_2)\right).
$$

Let us consider the cartesian diagram
\begin{equation*}
\xymatrix{
A^*(\M_0^{\leq 2}) \rat  \cart \ar[r]^{j^*} \ar[d]_{i^*} &\Q[\km_1, \km_2] \ar[d]^{q}\\
\Q[\km_1, \km_2, \gamma_2, \eta]/I \ar[r]^{p} & \Q[\km_1, \km_2, \eta]/\overline{I}
}
\end{equation*}
where $\overline{I}$ is the ideal generated by $\eta^2 +(2\km_2 +\km_1^2)\km_1^2$.

The map $q$ is injective so $i^*$ is injective too.

We set in $A^*(\M_0^{\leq 2}) \rat$ the classes $\gamma_2:=i_*1$ and $q:=i_*\eta$,
the ring we want (from injectivity of $i^*$) is isomorphic  to the subring of $A^*(\M_0^2)\rat$ generated by $\km_1, \km_2, \gamma_2, \gamma_2 \eta$ so we have
\begin{propos}
$$
A^*(\M_0^{\leq 2}) \rat = \frac{\Q [\km_1, \km_2, \gamma_2, q]}{(q^2 +\gamma_2 ^2(2\km_2 + \km_1^2)(\km_1^2-4\gamma_2))}
$$
\end{propos}

\newpage

\section {The third stratum}

The third stratum splits into two components.

\bigskip

{\bf The first component}
$\;$

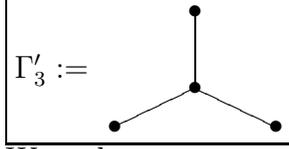
\framebox{$\Gamma'_3:=\;\vcenter{\small{\treb}}$}

We order components in $\Delta_1$. Let us note that the component corresponding to the central vertex has three points fixed by the other three components, consequently, given a permutation of the external vertices, there is an unique automorphism related to the central vertex.
 
The group $\text{Aut}(C^{\Gamma'_3})$ is therefore isomorphic to $S_3 \ltimes (E^3)$. We have
$$
A^*_{E^3} \rat \cong \Q [w_1, w_2, w_3],
$$
on which $S_3$ acts by permuting the three classes
\begin{eqnarray*}
w_1:=\psi^1_{\infty,1} \quad w_2:=\psi^1_{\infty,2} \quad w_3:=\psi^1_{\infty,3}
\end{eqnarray*}
So we have
$$
\begin{array}{l}
\phi^*\km_1=-(w_1 + w_2 + w_3),\\
\phi^*\km_2:=-(w_1^2 + w_2^2 + w_3^2),\\
\phi^*\km_3:=-(w_1^3 + w_2^3 + w_3^3);
\end{array}
$$
conesequently
$$
A^*(\M_0 ^{\Gamma'_3}) \rat \cong \Q [\km_1, \km_2, \km_3].
$$

We fix the following notation
\begin{equation*}
\xymatrix@=1.3pc{
\widetilde \M_0^{\Gamma'_3}   \ar@/_2pc/[dddrr]_\phi \ar[rr]^-{f}& &\Psi(\Gamma_2, \Gamma'_3)\ar[rr]^{pr_2} \ar [ddd]^{pr_1}& & \left( \widetilde{\M_0^{\Gamma_2}} \right) ^{\leq 2} \ar[ddd]^{\Pi}\\
& &\\
& & & \\
& & \M_0^{\Gamma'_3} \ar[rr]^{in}& &\M_0^{\leq 3}
}
\end{equation*}
where $f$ is the union of all $f_{\alpha}$.

First of all let us notice that the class $\phi^*\gamma_3':=\pi^*c_3(\Nf_{i})=w_1 w_2 w_3$
depends on Mumford classes in the following way $6\gamma'_3=-(\km_1^3 - 3 \km_1 \km_2 + 2 \km_3)$
so we can write
$$
A^*(\M_0^{\Gamma'_3}) \rat = \Q[\km_1, \km_2, \gamma'_3].
$$
The restriction of $\gamma_2$ to $A^*(\M_0^{\Gamma'_3})$is
\begin{eqnarray*}
\gamma_2:=pr_{1*}\frac{1}{2}c_2(pr_1^*(\Nf_{in})/\Nf_{pr_2})&=&\phi_*f^*\frac{1}{2}c_2(pr_1^*(\Nf_{in})/\Nf_{pr_2})\\
&=&\frac{1}{2}\pi_*(w_1 w_3)
\end{eqnarray*}
from which $\phi^*\gamma_2=w_1 w_3 + w_1 w_2 + w_2 w_3$ consequently, by writing $2\gamma_2 = \km_1^2 - \km_2$ we have
$$
A^*(\M_0^{\Gamma'_3}) \rat = \Q[\km_1, \gamma_2, \gamma'_3].
$$
In order to restrict the class $q$ let us notice that we can write
\begin{eqnarray*}
f^*r = 0 \quad f^*t_1= w_1 \quad f^*t_2 = w_3
\end{eqnarray*}
from which we have
$$
f^* \left( \frac{1}{2}(t_1 - t_2)(2t - t_1 + t_2) \right)= - \frac{1}{2}(w_1 - w_3)^2
$$
and so
\begin{eqnarray*}
\phi^*q&=&\phi^*\phi_* \left( - \frac{1}{2}(w_1 - w_3)^2 w_1w_3 \right)\\
&=& - ( (w_1 - w_3)^2 w_1 w_3 + (w_1 - w_2)^2 w_1 w_2 + (w_2 - w_3)^2 w_2 w_3)
\end{eqnarray*}
we can therefore write $q= - \gamma_2(\km_1^2 - 4 \gamma_2) + 3 \gamma'_3 \km_1$.
With reference to the inclusions
$$
\M_0^{\Gamma'_3} \xrightarrow{i} \M_0^{\leq 2} \cup \M_0^{\Gamma'_3} \xleftarrow{j} \M_0^{\leq 2}
$$
we have the fiber square:
\begin{equation*}
\xymatrix{
A^*(\M_0^{\leq 2} \cup \M_0^{\Gamma'_3}) \rat  \cart \ar[r]^{j^*} \ar[d]_{i^*} &\Q[\km_1, \km_2, \gamma_2, q]/I \ar[d]^{\varphi}\\
\Q[\km_1, \gamma_2, \gamma'_3] \ar[r]^{p} & \Q[\km_1, \gamma_2]
}
\end{equation*}
where  $I$ is the ideal generated by the polynomial $q^2 - \gamma_2^2(2\km_2 - \km_1^2)(\km_1^2 - 4\gamma_2)$ and the map $\varphi$ is surjective and such that
$$
\ker \varphi = (q + \gamma_2(\km_1^2 - 4 \gamma_2), \quad \km_2 + 2\gamma_2 - \km_1^2).
$$ 
Now let us observe that from Lemma (\ref{quoz}) the ring in question is isomorphic to
$$
A/(0,I):=\frac{\Q[\km_1, \gamma_2, \gamma'_3] \times_{\Q[\km_1, \gamma_2]} \Q[\km_1, \km_2, \gamma_2, q] }{(0,I)}.
$$
Set, with abuse of notation
\begin{eqnarray*}
\km_1:=(\km_1,\km_1) \quad & \gamma_2:=(\gamma_2,\gamma_2) \quad & \gamma'_3:=(\gamma'_3, 0)\\
\km_2:=(\km_1^2 - 2\gamma_2, \km_2) \quad & q:=(-\gamma_2(\km_1^2 - 4\gamma_2),q)&
\end{eqnarray*}
Straightforward arguments lead us to state the following
\begin{lemma}
The classes $\km_1,\gamma_2, \gamma'_3, \km_2, q$ generate the ring $A$.
\end{lemma}
Let us compute the ideal of relations. Let $T(\km_1, \gamma_2, \gamma'_3, \km_2, q)$ be a polynomial in $\Q[\km_1, \gamma_2, \gamma'_3, \km_2, q]$. It is zero in $A$ iff
\begin{eqnarray*}
&& T(\km_1, \gamma_2, \gamma'_3, \km_1^2 - 2\gamma_2, -\gamma_2(\km_1^2 - 4\gamma_2)) = 0 \text{ in } \Q[\km_1, \gamma_2, \gamma'_3]\\
&& T(\km_1, \gamma_2, 0, \km_2, q) = 0 \text{ in } \Q[\km_1, \gamma_2, \km_2, q]
\end{eqnarray*}

in particular this implies that $T$ is in the ideal of $\gamma'_3$. Consequently the polynomial $T=:\gamma'_3 \widehat{T}$ is zero in $A$ iff $\widehat{T}(\km_1, \gamma_2, \gamma'_3, \km_1^2 - 2\gamma_2, -\gamma_2(\km_1^2 - 4\gamma_2))$ is zero in  $\Q[\km_1, \gamma_2, \gamma'_3]$. The ideal of relations in $A$ is so generated by $\gamma'_3(-\km_2 + 2\gamma_2 - \km_1^2)$ and $\gamma'_3(q + \gamma_2(\km_1^2- 4\gamma_2))$.

Finally let us notice that the ideal $(0,I)$ is generated in $A$ by the polynomial $ q^2 + \gamma_2^2(2\km_2 + \km_1^2)(\km_1^2 - 4\gamma_2)$. So we can conclude that
$$
A^*(\M_0^{\leq 2} \cup \M_0^{\Gamma'_3}) \rat = \Q[\km_1, \km_2, \gamma_2, \gamma'_3, q]/J,
$$
where $J$ is the ideal generated by the polynomials
\begin{eqnarray*}
&& q^2 + \gamma_2^2(2\km_2 + \km_1^2)(\km_1^2 - 4\gamma_2)\\
&& \gamma'_3(- \km_2 + 2\gamma_2 - \km_1^2)\\
&&\gamma'_3(q + \gamma_2(\km_1^2- 4\gamma_2))
\end{eqnarray*}

{\bf The second component}

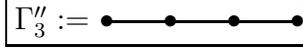
\framebox{$\Gamma''_3:=\;\vcenter{\small{\trea}}$}

The group of $\text{Aut}(C^{\Gamma'_3})$ is $\Cg_2 \ltimes (E \times E \times \Gm \times \Gm)$. The action of $\tau$ on this group exchange  simultaneously the components related to $\Gm$ and those related to $E$.

We have the isomorphism
$$
A^*_{ E \times E \times \Gm \times \Gm} \rat \cong \Q [v_1, \dots, v_4]
$$ 
where
\begin{eqnarray*}
v_1=\psi^1_{\infty,1} \quad v_2=\psi^1_{\infty, 2} \quad v_3=\psi^2_{\infty, 1} \quad v_4&=&\psi^2_{\infty, 2}
\end{eqnarray*}
by gluing curves such that the two central components corresponds in the point at infinity.

It follows that the action induced by $\tau$ is $$\tau(v_1,v_2,v_3,v_4)=(v_2, v_1,v_4,v_3)$$
Since $\Cg_2$ has order 2, the invariant polynomials are algebraically generated by the invariant polynomials of degree at most two (see. Theorem 7.5 \cite{CLO}). It is easy to see that a basis for the linear ones is given by $u_1:= v_1 + v_2, \; u_2:=v_3 + v_4$. For the vector subspace of invariant polynomials of degree two, we can compute a linear basis by using Reynolds' operator
$$
 \begin{array}{ll}
u_3:=v_1 ^2 +  v_2 ^ 2, & u_6:=v_1 v_3 + v_2 v_4,\\ 
u_4:=v_3 ^2 + v_4 ^2, & u_7:=v_1 v_2,\\
u_5:=v_1 v_4 + v_2 v_3, & u_8:=v_3 v_4
\end{array}
$$
now we note that
$$
u_6= u_1 u_2 - u_5,\quad u_7= (u_1^2 - u_3)/2,\quad u_8= (u_2^2 - u_4)/2
$$
Consequently we can write
$$
A^* \M_0^{\Gamma'_3} \rat \cong \Q [u_1, \dots, u_5]/I,
$$
where $I$ is the ideal generated by the polynomial
\begin{equation} \label{polgammauno}
2u_3 u_4 + 2u_1 u_2 u_5 -u_2 ^2 u_3 -u_1^2 u_4 - 2u_5^2
\end{equation}
With reference to the degree two covering  $\widetilde{\M}_0^{\Gamma'_3} \xrightarrow{\phi} \M_0^{\Gamma'_3}$
we have: $- u_1=\phi^*\km_1 \; \; -u_3=\phi^*\km_2$.
In order to compute the restriction of the closure of the classes $\gamma_2$ and $q$ of $\M_0^{\leq 2}$ let us fix the notation of the following diagram
\begin{equation*}
\xymatrix@=2.0pc{
& &\Ms_{0,1}^1 \times \Ms_{0,2}^0 \times \Ms_{0,1}^0 \ar[drr]^{pr_2^{'}} \\ \widetilde{\M}_0^{\Gamma'_3} \ar[urr]^{f_1} \ar[drr]^{f_3} \ar@/_2pc/[dddrr]_\phi \ar[rr]^-{f_2}& &(\Ms_{0,1}^0 \times \Ms_{0,2}^1 \times \Ms_{0,1}^0)^I \ar[rr]^{pr_2^{''}} & & \left( \widetilde{\M}_0^{\Gamma_2} \right)^{\leq 2} \ar[ddd]^{\Pi}\\
& &\Ms_{0,1}^0 \times \Ms_{0,2}^0 \times \Ms_{0,1}^1 \ar[dd]^{pr_1}\ar[urr]^{pr_2^{'''}}\\
& & & \\
& & \M_0^{\Gamma''_3} \ar[rr]^{i}& &\M_0^{\leq 3}
}
\end{equation*}
where $(\Ms_{0,1}^0 \times \Ms_{0,2}^1 \times \Ms_{0,1}^0)^I $ is the component where the marked points of the central curve (which is singular) are on different components.

As $\phi^*$ is an isomorphism to the algebra of polynomials which are invariants for the action of
$\Cg_2$, let us choose in $A^*(\M_0^{\Gamma''_3})\rat$ classes $\rho, \lambda, \mu$ such that
$u_2=\phi^*\rho \; \;,
u_4=\phi^*\lambda \; \;,
u_5=\phi^*\mu$.

First of all let us compute the restriction of the closure of $\gamma_2 \in A^*(\M_0^{\leq 2}) \rat$, the polynomial in the classes $\psi$ is $\frac{1}{2}$, we have
$c_3(\phi^*(\Nf_{i})) =(v_1-v_3)(v_3+v_4)(v_2-v_4)$, $c_1(f_1^*(\Nf_{pr_2^{'}}))=(v_1-v_3)$, $c_1(f_2^*(\Nf_{pr_2^{''}}))=(v_3+v_4)$ and $c_1(f_3^*(\Nf_{pr_2^{'''}}))=(v_2-v_4)$
from which we obtain the following relations
\begin{eqnarray*}
\phi^*\gamma'_3 &=&(v_1-v_3)(v_3+v_4)(v_2-v_4)\\
&=&\phi^*\left(\rho \left(\frac{1}{2}\left(\km_1^2 + \rho^2 - \km_2 - \lambda \right) -\mu \right) \right)\\
\phi^*\gamma_2&=&\left( (v_3+v_4)(v_2-v_4) + (v_1-v_3)(v_2-v_4) + (v_1-v_3)(v_3+v_4) \right)\\
&=&- \rho \km_1 - \mu + \frac{1}{2}(\km_1^2 - \rho^2 + \km_2 - \lambda)
\end{eqnarray*}
In order to have $\gamma_2$ among the generators, set
$$
\lambda= -2 \rho \km_1 - 2\mu + \km_1^2 - \rho^2 + \km_2 - 2 \gamma_2;
$$
the equation (\ref{polgammauno}) becomes
\begin{equation} \label{polbis}
K: \; \; \sigma^2-(2\km_2 + \km_1^2)((-\km_1 + 3\rho)(\km_1 + \rho) -4\gamma_2)
\end{equation}
where we have set $\sigma=2\mu - 2\km_2 - \km_1^2 + \km_1 \rho$.
Then we can write the ring $A^*(\M_0^{\Gamma''_3})\rat$ as $\Q[\km_1, \rho, \km_2, \gamma_2, \sigma]/K$.
In the new basis we have $\gamma'_3=\rho( \rho^2 - \rho \km_1 + \gamma_2)$.
Let us restricts the closure of $q \in A^*(\M_0^{\leq 2})\rat$ to $\M_0^{\Gamma_1}$. We recall that the related polynomial in classes $\psi$ of $\widetilde{\M}_0^{\Gamma_2}$ is $\frac{1}{2}(t_1-t_2)(2r - t_1+ t_2)$.
On $\Ms_{0,1}^1 \times \Ms_{0,2}^0 \times \Ms_{0,1}^0$ we have $
f_1^*t_1=v_3 \quad f_1^*t_2 = v_2 \quad f_1^*r =-v_4,$
on $(\Ms_{0,1}^0 \times \Ms_{0,2}^1 \times \Ms_{0,1}^0)^I$ we have $f_2^*t_1= v_1\quad f_2^*t_2 = v_2 \quad f_2^*r =\frac{\psi^2_\infty -\psi^2_0}{2}=\frac{v_3 - v_4}{2},$
and in the end on $\Ms_{0,1}^0 \times \Ms_{0,2}^0 \times \Ms_{0,1}^1$ we have $f_3^*t_1=v_1 \quad f_3^*t_2 =v_4 \quad f_3^*r =v_3,$
consequently the polynomials related to the three components of $\Psi(\Gamma_2, \Gamma''_3)$ are
$P_1 = \frac{1}{2}(v_3-v_2)(v_2 - v_3 - 2v_4)$, $P_2= \frac{1}{2}(v_1 - v_2) \left( v_4 - v_3 - v_1 + v_2 \right)$ and $P_3 = \frac{1}{2}(v_4 - v_1)(v_1 - 2v_3 - v_4)$
from which
\begin{eqnarray*} 
\phi^*q&=& 2\sum_{\alpha=1}^3 (P_\alpha N_\alpha)
= \phi^*((3\rho + \km_1)\gamma''_3 + (\rho^2 - \gamma_2)\sigma)
\end{eqnarray*}
this means that the image of $q$ in $B:=A^*(\M_0^{\Gamma''_3}) \rat / {\gamma''_3}$ is $(\rho^2 - \gamma_2)\sigma$.
Further let us notice that the image of $\gamma'_3$ in $B$ is 0.

In order to compute $A^*(\M_0^{\leq 3}) \rat$, let us consider its isomorphism with the fibered product
\begin{equation*}
\xymatrix{
 A^*(\M_0^{\leq 3}) \rat \cart \ar[r]^{j^*} \ar[d]_{i^*} & \Q [\km_1, \km_2, \gamma_2, \gamma'_3, q]/J \ar[d]^{\varphi}\\
\Q[\km_1, \rho, \km_2, \gamma_2, \sigma]/K \ar[r]^{p} & B
}
\end{equation*}
With reference to the fiber square
\begin{equation*}
\xymatrix{
 A \cart \ar[r]^{j^*} \ar[d]_{i^*} & \Q [\km_1, \km_2, \gamma_2, \gamma'_3, q] \ar[d]^{\varphi}\\
\Q[\km_1, \rho, \km_2, \gamma_2, \sigma]/K \ar[r]^{p} & B
}
\end{equation*}
we have that the ring in question is isomorphic to the quotient $A/(0,J)$.

Now we look for generators of the ring $A$.

Again a straightforward argument leads us to state
\begin{lemma}
The following elements of $A$
$$
\begin{array}{lll}
\km_1:=(\km_1, \km_1); & \km_2:=(\km_2,\km_2); & \gamma_2:=(\gamma_2);\\
\gamma'_3:=(0,\gamma'_3); & q:=((\rho^2 - \gamma_2)\sigma,q); & \gamma''_3:=(\gamma''_3,0);\\
r:=(\gamma''_3 \rho,0); & s:=(\gamma''_3 \sigma, 0) & t:=(\gamma''_3 \rho^2, 0);\\
& u:=(\gamma''_3 \rho \sigma,0)&
\end{array}
$$
are generators of the ring.
\end{lemma}

Now let us call $\widetilde A$ the ring $\Q[\km_1, \km_2, \gamma_2, \gamma'_3, \gamma''_3, q, r,s,t,u]$. We have by fact defined a surjective homomorphism $a: \widetilde A \to A$; we call again $j^*$ and $i^*$ their composition with $a$, we have $\ker a = \ker i^* \cap \ker j^*$.
Now, for what we've seen, we have:
\begin{eqnarray*}
\ker i^*&=& \left( 
	\begin{array}{l} 
	\gamma'_3;\\ 
	r^2 - \gamma''_3t; \\
	rs - \gamma''_3 u;\\
	s^2 - (\gamma''_3)^2(2\km_2 + \km_1^2)((-\km_1 + 2\rho)(\km_1 + \rho) + 4\gamma_2)
	\end{array}
	\right)\\
\ker j^* &=& (\gamma''_3, r,s,t,u)
\end{eqnarray*}
and so:
$$
\ker a= \left(
\begin{array}{l}
\gamma'_3 \gamma''_3; \quad \gamma'_3 r; \quad \gamma'_3 s; \quad \gamma'_3 t; \quad \gamma'_3 u;\\
r^2 - \gamma''_3t; \quad rs - \gamma''_3 u;\\
s^2 - (\gamma''_3)^2(2\km_2 + \km_1^2)((-\km_1 + 2\rho)(\km_1 + \rho) + 4\gamma_2)
\end{array}
\right)
$$

We make the quotient of $A=A'/ \ker a$ by the ideal $(0,J)$ and we obtain the following:
\begin{theorem}\label{finthr}
The ring $A^*(\M_0 ^{\leq 3}) \rat$ is: $$\Q[\km_1, \km_2, \gamma_2, \gamma'_3, \gamma''_3, q, r,s,t,u]/L$$
where $L$ is the ideal generated by the polynomials
$$
\begin{array}{l}
\gamma'_3(- \km_2 + 2\gamma_2 - \km_1^2); \quad \gamma'_3(q + \gamma_2(\km_1^2 - 4 \gamma_2));\\
\gamma'_3 \gamma''_3; \quad \gamma'_3 r; \quad \gamma'_3 s; \quad \gamma'_3 t; \quad \gamma'_3 u;\\
r^2 - \gamma''_3t; \quad rs - \gamma''_3 u;\\
q^2 + (\gamma_2)^2(2 \km_2 + \km_1^2)(\km_1^2 - 4 \gamma_2);\\
s^2 - (2\km_2 + \km_1^2)((- \km_1\gamma''_3 + 2r)(\km_1\gamma''_3 + r) + 4\gamma_2(\gamma''_3) ^2)
\end{array}
$$

\end{theorem}

\end{document}